\documentclass[11pt,a4paper,reqno]{amsart}
\usepackage{amsmath,amsfonts,amssymb,amsthm,mathrsfs,mathtools}
\usepackage{microtype,enumitem,comment,needspace}
\usepackage[dvipsnames]{xcolor}
\usepackage[pagebackref=true,colorlinks=true]{hyperref}
\hypersetup{linkcolor=blue,citecolor=blue,urlcolor=blue,filecolor=blue,pdfborder={0 0 0},pdftitle={Sharp matrix-weighted estimates for the Christ--Goldberg maximal operator},pdfauthor={Xi Cen and Pengfei Lu}}
\numberwithin{equation}{section}
\theoremstyle{plain}
\newtheorem{theorem}{Theorem}[section]
\newtheorem{problem}{Problem}
\newtheorem{lemma}[theorem]{Lemma}
\newtheorem{proposition}[theorem]{Proposition}
\newtheorem{corollary}[theorem]{Corollary}
\theoremstyle{definition}
\newtheorem{remark}[theorem]{Remark}
\newcommand{\R}{\mathbb R}
\newcommand{\C}{\mathbb C}
\newcommand{\one}{\mathbf 1}
\newcommand{\dd}{\,\mathrm d}
\newcommand{\norm}[1]{\lVert#1\rVert}
\newcommand{\op}{\mathrm{op}}
\DeclareRobustCommand{\Probref}[1]{\mbox{\hyperref[#1]{Problem~\ref*{#1}}}}
\DeclareRobustCommand{\Thmref}[1]{\mbox{\hyperref[#1]{Theorem~\ref*{#1}}}}
\DeclareRobustCommand{\Lemref}[1]{\mbox{\hyperref[#1]{Lemma~\ref*{#1}}}}
\DeclareRobustCommand{\Propref}[1]{\mbox{\hyperref[#1]{Proposition~\ref*{#1}}}}
\DeclareRobustCommand{\Corref}[1]{\mbox{\hyperref[#1]{Corollary~\ref*{#1}}}}
\DeclareRobustCommand{\Remref}[1]{\mbox{\hyperref[#1]{Remark~\ref*{#1}}}}
\DeclareRobustCommand{\Secref}[1]{\mbox{\hyperref[#1]{Section~\ref*{#1}}}}
\renewcommand{\backref}[1]{}

\renewcommand{\backrefalt}[4]{\ifnum#1>0\relax\space{\footnotesize\color{blue}#2}\fi}
\renewcommand*{\eqref}[1]{\mbox{\hyperref[#1]{(\ref*{#1})}}}

\allowdisplaybreaks[1]
\title[Sharp matrix-weighted maximal estimates]{\bf Sharp matrix-weighted estimates for the Christ--Goldberg maximal operator}
\author[X. Cen]{Xi Cen}
\address{Xi Cen\\
School of Science\\
China University of Mining and Technology-Beijing\\
Beijing 100083, People's Republic of China}
\email{xicenmath@gmail.com}
\author[P. Lu]{Pengfei Lu}
\address{Pengfei Lu\\
Department of Mathematics\\
Northeastern University\\
Shenyang 110004, People's Republic of China}
\email{pflumath@gmail.com}
\date{}
\begin{document}

\begin{abstract}
We resolve the sharpness conjecture of Lerner, Li, Ombrosi and Rivera-R\'ios~\cite[\emph{Ann. Sc. Norm. Super. Pisa Cl. Sci.}, 2024, p.~3]{LLOR} for the Christ--Goldberg maximal operator when $2\le p<\infty$. Under the constraint $[W]_{A_p}\le t$, its optimal multiplier weak-type norm has order $t^{1/(p-1)}$, with no logarithmic loss. The matching lower bound holds even on the standard dyadic grid for smooth positive periodic scalar weights. For every $1<p<\infty$, we further determine the optimal strong- and weak-type norms when the averaging family has overlap at most $K$. Their orders are $\min\{(Kt)^{1/p},t^{1/(p-1)}\}$ and $\min\{(Kt)^{1/p},t^{1/(p-1)},t^{2/p}\}$, respectively. For each $t\ge1$, the same scalar weight and unit function realize both orders, up to constants, on increasing dyadic chains for all $K\ge1$.
\end{abstract}
\maketitle

\section{Introduction}\label{sec:intro}
\subsection{Background and motivations}\label{subsec:background}
Quantitative weighted inequalities ask how the norm of an operator depends on the characteristic of a weight. For the Christ--Goldberg maximal operator, the sharp strong-type estimate is known, but a logarithmic gap remained in the multiplier weak-type estimate for $p\ge2$~\cite[(1.4) and Theorem~1.3]{LLOR}. We remove this gap. We also determine how the optimal strong and weak norms change when the number of averaging cubes containing a point is bounded.

Let $n\in\mathbb N=\{1,2,\ldots\}$ and $1<p<\infty$. For $f\in L^1_{\mathrm{loc}}(\R^n)$ the uncentered Hardy--Littlewood maximal operator is $Mf(x)=\sup_{Q\ni x}\langle |f|\rangle_Q$, where $Q$ ranges over axis-parallel cubes and $\langle f\rangle_Q=|Q|^{-1}\int_Q f$. Here $|E|$ is Lebesgue measure. For a locally integrable weight $w>0$ almost everywhere, set $\norm f_{L^p(w)}=(\int_{\R^n}|f|^pw\dd x)^{1/p}$ and $p'=p/(p-1)$. Muckenhoupt~\cite{Muck} characterized the boundedness of $M$ on $L^p(w)$ by
\begin{equation}\label{eq:scalarAp}
 [w]_{A_p}=\sup_Q\langle w\rangle_Q\langle\sigma_w\rangle_Q^{p-1}<\infty,\qquad \sigma_w=w^{-1/(p-1)}.
\end{equation}
We write $w\in A_p$ when \eqref{eq:scalarAp} holds. Buckley~\cite[Theorem~2.5]{Buckley} proved the strong bound of order $[w]_{A_p}^{1/(p-1)}$ with optimal exponent. For several weights, Li, Moen and Sun~\cite{LMS} obtained sharp estimates for multilinear maximal functions, while Cao and Xue~\cite{CaoXue} studied the two-weight characterization of multilinear fractional maximal operators.

In the scalar case, conjugating $M$ by the isometry $g\mapsto w^{1/p}g$ from $L^p(w)$ onto $L^p=L^p(\R^n)$ leads to the multiplier maximal operator
\begin{equation}\label{eq:scalarconj}
 M_{w,p}f(x)=w(x)^{1/p}M(w^{-1/p}f)(x).
\end{equation}
For $w\in A_p$, $\norm{M}_{L^p(w)\to L^p(w)}=\norm{M_{w,p}}_{L^p\to L^p}$. The corresponding weak norms need not agree. With $w(E)=\int_Ew\dd x$, the ordinary and multiplier weak quantities for $g\in L^p(w)$ are
\[
 \sup_{\lambda>0}\lambda w(\{Mg>\lambda\})^{1/p},\qquad
 \sup_{\lambda>0}\lambda|\{w^{1/p}Mg>\lambda\}|^{1/p}.
\]
By \eqref{eq:scalarconj}, the second is $\norm{M_{w,p}(w^{1/p}g)}_{L^{p,\infty}}$, with $\norm h_{L^{p,\infty}}=\sup_{\lambda>0}\lambda|\{x:|h(x)|>\lambda\}|^{1/p}$. Muckenhoupt and Wheeden~\cite{MW} studied multiplier weak inequalities on the line. Sawyer~\cite{Sawyer} introduced mixed weak inequalities, further developed by Cruz-Uribe, Martell and P\'erez~\cite{CMP} and Li, Ombrosi and P\'erez~\cite{LOP}. In the latter setting the target measure is weighted as well. Our target space $L^{p,\infty}=L^{p,\infty}(\R^n)$ uses Lebesgue measure. This formulation also accommodates matrix weights.

Christ and Goldberg~\cite{CG} introduced a maximal operator adapted to matrix $A_2$ weights, and Goldberg~\cite{Goldberg} extended the theory to matrix $A_p$ weights for $1<p<\infty$. Let $W$ be a locally integrable Hermitian $d\times d$ matrix weight, $d\ge1$, positive definite almost everywhere. Using spectral powers and the Euclidean operator norm $\norm{\cdot}_{\op}$, set
\begin{equation}\label{eq:matrixAp}
 [W]_{A_p}=\sup_Q\frac1{|Q|}\int_Q\left(\frac1{|Q|}\int_Q
 \norm{W(x)^{1/p}W(y)^{-1/p}}_{\op}^{p'}\dd y\right)^{p/p'}\dd x.
\end{equation}
The expression for a fixed cube $Q$ in \eqref{eq:matrixAp} is denoted by $[W]_{A_p,Q}$.
For $W\in A_p$, meaning $[W]_{A_p}<\infty$, and $\mathbf f\in L^p(\R^n,\C^d)$, the Christ--Goldberg maximal operator, in the notation of~\cite[p.~2]{LLOR}, is defined by
\begin{equation}\label{eq:matrixconj}
 M_{W,p}\mathbf f(x)=\sup_{Q\ni x}\frac1{|Q|}\int_Q
 |W(x)^{1/p}W(y)^{-1/p}\mathbf f(y)|\dd y.
\end{equation}
The domain consists of $\C^d$-valued functions, while $M_{W,p}\mathbf f$ is nonnegative and scalar-valued. Here $|\cdot|$ is the Euclidean norm. If $w\in A_p$ is scalar and $I_d$ denotes the $d\times d$ identity matrix, then
\begin{equation}\label{eq:conj}
 [wI_d]_{A_p}=[w]_{A_p},\qquad M_{wI_d,p}\mathbf f=M_{w,p}(|\mathbf f|).
\end{equation}
Scalar lower bounds therefore transfer to every matrix dimension by using $f\mathbf e$ with $|\mathbf e|=1$. For an operator $T:X\to Y$ we write $\norm T_{X\to Y}=\sup_{\norm f_X\le1}\norm{Tf}_Y$, also when $Y$ is equipped with a weak quasi-norm. The symbols $\lesssim_\Theta$ and $\eqsim_\Theta$ denote upper and two-sided estimates with constants depending only on $\Theta$.

Isralowitz and Moen~\cite[Theorem~1.3, $\alpha=0$, $q=p$]{IM} established the strong estimate for the matrix operator in dimension $n\ge2$ (see also~\cite[Remark~3.5]{LLOR}). With the reduction in \Lemref{lem:dim}, it reads
\begin{equation}\label{eq:knownstrong}
 \norm{M_{W,p}}_{L^p(\R^n,\C^d)\to L^p}\lesssim_{n,d,p}[W]_{A_p}^{1/(p-1)}\qquad(1<p<\infty).
\end{equation}
At the endpoint $p=1$, Cruz-Uribe, Isralowitz, Moen, Pott and Rivera-R\'ios~\cite{CruzUribeEtAl2021} obtained a weak bound quadratic in their matrix $A_1$ characteristic. Lerner, Li, Ombrosi and Rivera-R\'ios~\cite{LLORendpoint} proved sharpness using scalar weights. For $p>1$, Cruz-Uribe and Sweeting~\cite{CS} obtained a weak bound with power $1+1/p$. Lerner, Li, Ombrosi and Rivera-R\'ios~\cite[Theorem~1.1]{LLOR} lowered the power to the sharp value $2/p$ when $1<p<2$. Since $2/p<1/(p-1)$ precisely in this range, it is the remaining case $p\ge2$ that calls for a sharpness test of \eqref{eq:knownstrong} in the weak norm.

\subsection{The original problem and the main results}\label{subsec:problem}
Lerner, Li, Ombrosi and Rivera-R\'ios~\cite[(1.4), p.~2]{LLOR} deduced from the strong bound that
\begin{equation}\label{eq:originalbound}
 \norm{M_{W,p}}_{L^p(\R^n,\C^d)\to L^{p,\infty}}\lesssim_{n,d,p}[W]_{A_p}^{1/(p-1)}\quad(1<p<\infty).
\end{equation}
Immediately before Theorem~1.3 on p.~3 of~\cite{LLOR}, they conjectured that their estimate (1.4) is sharp for $p\ge2$. In their notation, the one-dimensional scalar version of the problem reads as follows.
\begin{problem}[{Lerner--Li--Ombrosi--Rivera-R\'ios~\cite[p.~3, before Theorem~1.3]{LLOR}}]\label{prob:orig}
Let $2\le p<\infty$. Is it true that
\begin{equation}\label{eq:originalquestion}
 \varphi_{M_p}(t):=\sup_{[w]_{A_p(\R)}\le t}\norm{M_{w,p}}_{L^p(\R)\to L^{p,\infty}(\R)}\gtrsim_p t^{1/(p-1)}\quad(t\ge1)?
\end{equation}
\end{problem}
By~\cite[Theorem~1.3]{LLOR} and \eqref{eq:originalbound},
\begin{equation}\label{eq:prior-gap}
 \frac{t^{1/(p-1)}}{(\log(e+t))^{1/p}}\lesssim_p\varphi_{M_p}(t)\lesssim_p t^{1/(p-1)}.
\end{equation}
The lower bound in \eqref{eq:prior-gap} already excludes every smaller power of $t$. The remaining question is whether the factor $(\log(e+t))^{-1/p}$ can be removed.

Let $\mathcal D=\{2^j(k+[0,1)^n):j\in\mathbb Z,\ k\in\mathbb Z^n\}$ be the standard dyadic grid. For a countable family $\mathcal F$ of cubes, we use the same maximal operator restricted to $\mathcal F$:
\[
 M_{W,p}^{\mathcal F}\mathbf f(x)=\sup_{\substack{Q\in\mathcal F\\Q\ni x}}\frac1{|Q|}\int_Q
 |W(x)^{1/p}W(y)^{-1/p}\mathbf f(y)|\dd y,
\]
with the empty supremum equal to zero. We use $M_{w,p}^{\mathcal F}$ for its scalar counterpart. The weight characteristics still use all cubes. Write $C_c^\infty(U)$ for smooth functions compactly supported in $U$.

\begin{theorem}[Sharp matrix-weighted weak-type estimate]\label{thm:scalar}
Let $n,d\in\mathbb N$, $2\le p<\infty$ and $t\ge1$. Then
\begin{equation}\label{eq:matrixmain}
 \sup_{[W]_{A_p}\le t}\norm{M_{W,p}}_{L^p(\R^n,\C^d)\to L^{p,\infty}}\eqsim_{n,d,p}t^{1/(p-1)}.
\end{equation}
More precisely, even the scalar operator restricted to the standard dyadic grid satisfies
\begin{equation}\label{eq:main}
 \sup_{[w]_{A_p}\le t}\norm{M_{w,p}^{\mathcal D}}_{L^p\to L^{p,\infty}}\eqsim_{n,p}t^{1/(p-1)}.
\end{equation}
For every $t\ge1$, there exist a smooth positive weight $w_t(x)=u_t(x_1)$ with $[w_t]_{A_p}\le t$ and $u_t$ $1$-periodic, and a nonnegative $L^p$ unit function $f_t\in C_c^\infty((0,1/2)^n)$ such that
\[
 \norm{M_{w_t,p}^{\mathcal D}f_t}_{L^{p,\infty}}\gtrsim_{n,p}t^{1/(p-1)}.
\]
\end{theorem}

For $n=1$, \eqref{eq:main} gives the lower bound \eqref{eq:originalquestion} and hence answers \Probref{prob:orig}. Together with the scalar embedding \eqref{eq:conj} and~\cite[Theorem~1.1]{LLOR} for $1<p<2$, it gives, for every $n,d\ge1$, $1<p<\infty$ and $t\ge1$,
\begin{equation}\label{eq:allp}
 \sup_{[W]_{A_p}\le t}\norm{M_{W,p}}_{L^p(\R^n,\C^d)\to L^{p,\infty}}\eqsim_{n,d,p}
 t^{\min\{2/p,1/(p-1)\}}.
\end{equation}
In particular, for every matrix $A_p$ weight $W$ and $1<p<\infty$,
\[
 \norm{M_{W,p}}_{L^p(\R^n,\C^d)\to L^{p,\infty}}\lesssim_{n,d,p}[W]_{A_p}^{\min\{2/p,1/(p-1)\}}.
\]
For $1<p<2$, the exact order was established in~\cite[Theorem~1.1]{LLOR}. For $p\ge2$, \eqref{eq:knownstrong} supplies the upper bound, and \Thmref{thm:scalar} supplies the matching lower bound. In particular, \eqref{eq:allp} determines the full optimal growth, uniformly over the matrix weights with characteristic at most $t$.

Write $\one_E$ for the indicator of a measurable set $E$. For a countable family $\mathcal F$, its overlap is bounded by $K$ if
\begin{equation}\label{eq:depth}
 \left\|\sum_{Q\in\mathcal F}\one_Q\right\|_{L^\infty}\le K.
\end{equation}
For $\mathcal F=\{Q\}$, H\"older's inequality and its equality case give
\[
 \norm{M_{w,p}^{\{Q\}}}_{L^p\to L^p}=\langle w\rangle_Q^{1/p}\langle\sigma_w\rangle_Q^{1/p'}\le [w]_{A_p}^{1/p}.
\]
Summing the $p$th powers over cubes satisfying \eqref{eq:depth} gives the bound $(Kt)^{1/p}$ under $[w]_{A_p}\le t$. The unrestricted operator gives another bound, independent of $K$. We determine the exact dependence on both parameters and show that the two optimal orders can be realized simultaneously.

\begin{theorem}[Sharp matrix-weighted estimates for restricted maximal operators]\label{thm:depth}
Let $n,d\ge1$, $1<p<\infty$, $t\ge1$ and $K\in\mathbb N$. The suprema below range over all matrix $A_p$ weights and countable cube families. Then
\begin{align}
 \sup_{\substack{[W]_{A_p}\le t\\\|\sum_{Q\in\mathcal F}\one_Q\|_{L^\infty}\le K}}\norm{M_{W,p}^{\mathcal F}}_{L^p(\R^n,\C^d)\to L^p}
 &\eqsim_{n,d,p}\min\{(Kt)^{1/p},t^{1/(p-1)}\},\label{eq:depthstrong}\\
 \sup_{\substack{[W]_{A_p}\le t\\\|\sum_{Q\in\mathcal F}\one_Q\|_{L^\infty}\le K}}\norm{M_{W,p}^{\mathcal F}}_{L^p(\R^n,\C^d)\to L^{p,\infty}}
 &\eqsim_{n,d,p}\min\{(Kt)^{1/p},t^{1/(p-1)},t^{2/p}\}.\label{eq:depthsharp}
\end{align}
For each fixed $n,p,t$, there exist a scalar weight $w_t$ with $[w_t]_{A_p}\le t$, a nonnegative $L^p$ unit function $h_t$, and finite dyadic chains $\mathcal F_K\subseteq\mathcal F_{K+1}$, $\#\mathcal F_K\le K$, realizing both lower bounds on $h_t$ uniformly in $K,t$. The union of the chains is finite. For matrix weights the same lower bounds hold for $W=w_tI_d$ and $h_t\mathbf e$, where $|\mathbf e|=1$.
\end{theorem}

\Needspace{18\baselineskip}
These suprema yield the following sharp inequalities.
\begin{corollary}[Sharp inequalities and separation of the norms]\label{cor:ratio}
Let $n,d\ge1$, $1<p<\infty$, $W\in A_p$ and let $\mathcal F$ be a countable family of cubes satisfying \eqref{eq:depth} for some $K\in\mathbb N$. Then
\begin{align}
 \norm{M_{W,p}^{\mathcal F}}_{L^p(\R^n,\C^d)\to L^p}
 &\lesssim_{n,d,p}\min\{(K[W]_{A_p})^{1/p},[W]_{A_p}^{1/(p-1)}\},\label{eq:actualstrong}\\
 \norm{M_{W,p}^{\mathcal F}}_{L^p(\R^n,\C^d)\to L^{p,\infty}}
 &\lesssim_{n,d,p}\min\{(K[W]_{A_p})^{1/p},[W]_{A_p}^{1/(p-1)},[W]_{A_p}^{2/p}\}.\label{eq:actualweak}
\end{align}
On the examples of \Thmref{thm:depth}, $[w_t]_{A_p}\eqsim_{n,p}t$, and both upper bounds are attained up to constants, also by the same unit function $h_t$. Moreover,
\begin{equation}\label{eq:ratio}
 \frac{\norm{M_{w_t,p}^{\mathcal F_K}}_{L^p\to L^p}}{\norm{M_{w_t,p}^{\mathcal F_K}}_{L^p\to L^{p,\infty}}}
 \eqsim_{n,p}\begin{cases}
 \left(\dfrac{\min\{K,t^{1/(p-1)}\}}{\min\{K,t\}}\right)^{1/p},&1<p<2,\\[3pt]
 1,&p\ge2.
 \end{cases}
\end{equation}
The same equivalence holds for the ratio of the norms evaluated at $h_t$. In particular, when $1<p<2$ and $K\ge t^{1/(p-1)}$, both ratios are comparable to $t^{(2-p)/(p(p-1))}$.
\end{corollary}

\begin{remark}[Sharpness and saturation]\label{rem:mechanism}
For the smooth examples of \Thmref{thm:scalar}, the lower bound in \eqref{eq:main} and \eqref{eq:originalbound} imply $[w_t]_{A_p}\eqsim_{n,p}t$. Hence, for $p\ge2$,
\begin{equation}\label{eq:smoothactual}
 \norm{M_{w_t,p}^{\mathcal D}}_{L^p\to L^{p,\infty}}\eqsim_{n,p}[w_t]_{A_p}^{1/(p-1)}.
\end{equation}
Thus \eqref{eq:originalbound} admits no uniform improvement by a factor tending to zero as $[W]_{A_p}\to\infty$, even for smooth scalar weights. For any prescribed $K$, the overlap can also be made exactly $K$ by adjoining $K-\#\mathcal F_K$ distinct dyadic ancestors of the largest cube of the witnessing chain. The lower bounds persist, and the upper bounds apply with overlap exactly $K$. Comparing the terms in \eqref{eq:depthstrong} and \eqref{eq:depthsharp}, the strong norm reaches its unrestricted order at $K\eqsim t^{1/(p-1)}$, and the weak norm does so at $K\eqsim\min\{t,t^{1/(p-1)}\}$.

The distinction between strong and weak norms is reflected in the lower-bound constructions. For $p\ge2$, \Propref{prop:Ap} gives $N$ disjoint intervals of length $N$ on which the maximal function of a unit function is bounded below by $c_pN^{(p-2)/p}$ when $B=2$. The corresponding weak norm is at least $c_pN$, while the characteristic has order $N^{p-1}$. For $1<p<2$, we retain the intervals used in~\cite[Theorem~1.1]{LLOR} for the weak norm and modify the weight near the origin to obtain the additional strong norm growth.
\end{remark}

\subsection{Organization}\label{subsec:organization}
\Secref{sec:mainproof} proves the main results using the auxiliary lemmas stated at its beginning. Their proofs, including the full characteristic estimates and the passage from intervals to cubes, appear in \Secref{sec:aux}.

\subsection{Notation}\label{subsec:notation}
The unweighted norm on $L^p(\R^n,\C^d)$ is $\norm{\mathbf f}_{L^p(\R^n,\C^d)}=(\int_{\R^n}|\mathbf f|^p\dd x)^{1/p}$, with the analogous scalar norm on $L^p$. On $\R_+=(0,\infty)$, the $A_p$ characteristic in \eqref{eq:scalarAp} uses intervals $(a,b)$ with $0\le a<b<\infty$. Restricted scalar operators on $\R_+$ use the same interval averages. For a single cube $Q$, the restricted maximal operators reduce to the averages
\begin{equation}\label{eq:scalarop}
 \begin{aligned}
 M_{v,p}^{\{Q\}}f(x)&=\one_Q(x)v(x)^{1/p}\langle v^{-1/p}|f|\rangle_Q,\\
 M_{W,p}^{\{Q\}}\mathbf f(x)&=\frac{\one_Q(x)}{|Q|}\int_Q|W(x)^{1/p}W(y)^{-1/p}\mathbf f(y)|\dd y.
 \end{aligned}
\end{equation}
In particular, $M_{W,p}^{\mathcal F}\mathbf f=\sup_{Q\in\mathcal F}M_{W,p}^{\{Q\}}\mathbf f$ pointwise, and likewise in the scalar case. For $W\in A_p$, choose $x_0\in Q$ with $W(x_0)$ invertible and the inner integral in \eqref{eq:matrixAp} finite. Then
\[
 \norm{W(y)^{-1/p}}_{\op}\le\norm{W(x_0)^{-1/p}}_{\op}\norm{W(x_0)^{1/p}W(y)^{-1/p}}_{\op}
\]
shows $\norm{W^{-1/p}}_{\op}\in L^{p'}(Q)$. H\"older gives local integrability of $W^{-1/p}\mathbf f$. For a unit vector $\mathbf e$, the scalar weight $a(x)=|W(x)^{1/p}\mathbf e|^p$ satisfies $1\le[a]_{A_p}\le[W]_{A_p}$, since $a(x)^{1/p}/a(y)^{1/p}\le\norm{W(x)^{1/p}W(y)^{-1/p}}_{\op}$. Thus the averages exist independently of representatives. Approximation by enclosing rational cubes gives measurability, and \eqref{eq:knownstrong} gives almost everywhere finiteness. These facts include $d=1$. Norms on intervals and half-lines use the indicated domains. We write $A\gtrsim_\Theta B$ for $B\lesssim_\Theta A$. Constants are uniform in the weights and scales.

Write $\R_+^n=(0,\infty)^n$, $q(x)=|x|_\infty^n$, where $|x|_\infty=\max_i|x_i|$, and $I(r)=(0,r)$, $Q(r)=[0,r^{1/n})^n$ for $r>0$. We write $\operatorname{dist}(x,E)=\inf_{y\in E}|x-y|$. For a measurable set $C\subseteq\R_+$ and a scalar weight $v$, write $\sigma_v=v^{-1/(p-1)}$ and $S_C=\int_C\sigma_v$. If $0<S_C<\infty$, set $h_C=S_C^{-1/p}\sigma_v^{1/p}\one_C$. All logarithms are natural, and $\log_B r=\log r/\log B$.

\section{Proofs of the main results}\label{sec:mainproof}
\subsection{Auxiliary lemmas}\label{subsec:auxstatements}
The auxiliary results below separate the estimates on the line, their transfer to cubes, and the modification of weights. Their proofs are given in \Secref{sec:aux}.

\begin{lemma}[A common test on an interval chain]\label{lem:common}
Let $v\in A_p(\R_+)$, $1<p<\infty$, $0=r_0<r_1<\cdots<r_J$, and $\mathcal I=\{I(r_j):1\le j\le J\}$. Suppose that $C\subseteq I(r_1)$ with $0<S=S_C<\infty$. Then $\norm{h_C}_{L^p(\R_+)}=1$. If an interval $I$ contains $C$ up to a null set, then
\begin{equation}\label{eq:sourceaverage}
 M_{v,p}^{\{I\}}h_C(x)=\frac{S^{1/p'}v(x)^{1/p}}{|I|}\quad(x\in I).
\end{equation}
Consequently, with $r_{j-1}<x<r_j$,
\begin{equation}\label{eq:chainidentity}
 \begin{aligned}
 M_{v,p}^{\mathcal I}h_C(x)&=\frac{S^{1/p'}v(x)^{1/p}}{r_j},\\
 \norm{M_{v,p}^{\mathcal I}h_C}_{L^p}^p&=S^{p-1}\sum_{j=1}^Jr_j^{-p}\int_{r_{j-1}}^{r_j}v.
 \end{aligned}
\end{equation}
\begin{equation}\label{eq:weakchain}
 \norm{M_{v,p}^{\mathcal I}h_C}_{L^{p,\infty}}^p
 = S^{p-1}\sup_{a>0}a^p\sum_{j=1}^J
 |\{x\in(r_{j-1},r_j):v(x)>(ar_j)^p\}|.
\end{equation}
For disjoint measurable $E_j\subseteq I(r_j)$ and $a>0$, it follows that
\begin{align}
 \norm{M_{v,p}^{\mathcal I}h_C}_{L^p}^p&\ge S^{p-1}\sum_{j=1}^Jr_j^{-p}\int_{E_j}v,\label{eq:strongsource}\\
 \norm{M_{v,p}^{\mathcal I}h_C}_{L^{p,\infty}}^p&\ge a^pS^{p-1}\sum_{j=1}^J|E_j|
 \quad\text{if }v^{1/p}\ge ar_j\text{ a.e. on }E_j.\label{eq:chainnorms}
\end{align}
All pointwise identities are understood almost everywhere.
\end{lemma}

\begin{lemma}[Reduction to dimension one]\label{lem:dim}
Let $1<p<\infty$ and $W\in A_p(\R)$ be a $d\times d$ matrix weight. Suppose that \eqref{eq:knownstrong} holds in dimension two with constant $C_{2,d,p}$. Then
\[
 \norm{M_{W,p}}_{L^p(\R,\C^d)\to L^p(\R)}\le C_{2,d,p}[W]_{A_p(\R)}^{1/(p-1)}.
\]
\end{lemma}

\begin{lemma}[Transfer of interval chains to dyadic cubes]\label{lem:radial}
Let $1<p<\infty$, $v\in A_p(\R_+)$, $B=2^n$ and $V(x)=v(q(x))$ away from zero, with $V(0)=1$. Then
\[
 [V]_{A_p(\R^n)}\lesssim_{n,p}[v]_{A_p(\R_+)}.
\]
Let $f\in L^p(\R_+)$ and define $F(x)=f(q(x))$ on $\R_+^n$ and zero elsewhere. Let $\mathcal I=\{I(r_j)\}_j$ with $r_j\in B^{\mathbb Z}$ and $\mathcal F=\{Q(r_j)\}_j$. Then
\begin{equation}\label{eq:transfer}
 \norm F_{L^p(\R^n)}=\norm f_{L^p(\R_+)},\qquad
 \norm{M_{V,p}^{\mathcal F}F}_{Y(\R^n)}=\norm{M_{v,p}^{\mathcal I}f}_{Y(\R_+)}
 \quad(Y=L^p,\ L^{p,\infty}).
\end{equation}
For finite chains, the overlap equals the cardinality.
\end{lemma}

\begin{proposition}[The periodic weight]\label{prop:Ap}
Let $p\ge2$, $N\in\mathbb N$ and $B\ge2$. Define
\[
 b=B^{\lceil\log_B(256BN)\rceil},\qquad s_k=bB^k\ (0\le k<2N),\qquad L=bB^{2N}.
\]
Define the disjoint intervals $C_k$ and $E_k$ by
\begin{equation}\label{eq:CE}
 \begin{aligned}
 C_k&=[5s_k/4,3s_k/2]\ (0\le k<N),& C&=\bigcup_{k=0}^{N-1}C_k,\\
 E_k&=[s_k,s_k+N]\ (0\le k<2N),& E&=\bigcup_{k=N}^{2N-1}E_k.
 \end{aligned}
\end{equation}
With $\Gamma=[0,b]\cup\bigcup_{k=0}^{2N-1}E_k\cup\{L\}$ and $r(z)=\max\{b,z\}$, set
\begin{equation}\label{eq:weight}
 \begin{aligned}
 w_0(z)&=Nr(z)^{p-1}+\frac{r(z)^p}{1+\operatorname{dist}(z,\Gamma)}\quad(0\le z\le L),\\
 w(z)&=w_0(\operatorname{dist}(z,2L\mathbb Z))\quad(z\in\R).
 \end{aligned}
\end{equation}
Then $w$ is continuous, positive and $2L$-periodic. Its characteristic satisfies
\begin{equation}\label{eq:Apglobal}
 [w]_{A_p(\R)}\eqsim_{B,p}N^{p-1}.
\end{equation}
Moreover, writing $S_C=\int_C w^{-1/(p-1)}$, we have
\begin{equation}\label{eq:testnorm}
 S_C\eqsim_{B,p}N^{(p-2)/(p-1)},\qquad w(z)\ge s_j^p\ (z\in E_j,\ N\le j<2N),\quad |E|=N^2.
\end{equation}
\end{proposition}

\begin{lemma}[Smooth periodic realization]\label{lem:smooth}
Let $N\in\mathbb N$, $p\ge2$ and $n\ge1$. Then there are a weight $w_*(x)=u_*(x_1)$ with $u_*$ smooth, positive and $1$-periodic, and a nonnegative $f_*\in C_c^\infty((0,1/2)^n)$ with $\|f_*\|_{L^p}=1$ such that
\[
 [w_*]_{A_p}\le K_{n,p}N^{p-1},\qquad
 \|M_{w_*,p}^{\mathcal D}f_*\|_{L^{p,\infty}}\gtrsim_{n,p}N.
\]
The constant $K_{n,p}\ge1$ depends only on $n$ and $p$.
\end{lemma}

For $1<p<2$ and integers $m\ge100$, let $u_m$ be the periodic scalar weight in~\cite[Section~3, pp.~10--13]{LLOR}, with their parameter $N=m$. Their construction on p.~11 and the bound (3.4) on p.~12 give
\begin{equation}\label{eq:oldprofile}
 [u_m]_{A_p(\R)}\lesssim_p m,\qquad u_m=1\ \text{on }[0,1],
\end{equation}
and
\begin{equation}\label{eq:oldplateaux}
 u_m=2^{(k+1)p}\quad\text{on }I_k^*=[2^k,2^k+k),\qquad 3\le k\le m .
\end{equation}
\begin{lemma}[Modification near the origin]\label{lem:patch}
Let $1<p<2$, $m\in\mathbb N$, $m\ge100$, and let $u_m$ have the properties \eqref{eq:oldprofile} and \eqref{eq:oldplateaux}. Write $\varepsilon=m^{-1/(p-1)}$ and $\beta=(p-1)(1-\varepsilon)$. Let $0<\tau\le\varepsilon$ and define
\[
 v(z)=\begin{cases}(z/\tau)^\beta,&0<z<\tau,\\u_m(z),&z\ge\tau.\end{cases}
\]
Then $[v]_{A_p(\R_+)}\lesssim_p m$.
\end{lemma}

\subsection{Upper bounds}\label{subsec:upper}
By Chebyshev's inequality and \eqref{eq:knownstrong},
\[
 \begin{aligned}
 \norm{M_{W,p}^{\mathcal D}}_{L^p(\R^n,\C^d)\to L^{p,\infty}}
 &\le\norm{M_{W,p}}_{L^p(\R^n,\C^d)\to L^{p,\infty}}\\
 &\le\norm{M_{W,p}}_{L^p(\R^n,\C^d)\to L^p}
 \lesssim_{n,d,p}[W]_{A_p}^{1/(p-1)}.
 \end{aligned}
\]
Taking the suprema gives both upper bounds in \Thmref{thm:scalar}.

\begin{proof}[Upper bounds in \Thmref{thm:depth}]
Fix $W\in A_p$, a countable $\mathcal F$ satisfying \eqref{eq:depth}, and $\mathbf h\in L^p(\R^n,\C^d)$. For almost every $x\in Q$, H\"older's inequality gives
\[
 \begin{aligned}
 |M_{W,p}^{\{Q\}}\mathbf h(x)|^p
 &\le\left\langle\norm{W(x)^{1/p}W(\cdot)^{-1/p}}_{\op}^{p'}\right\rangle_Q^{p/p'}
 \left\langle|\mathbf h|^p\right\rangle_Q.
 \end{aligned}
\]
Integrating over $Q$ and using Tonelli's theorem, we obtain
\[
 \begin{aligned}
 \int_Q|M_{W,p}^{\{Q\}}\mathbf h|^p
 &\le [W]_{A_p,Q}\int_Q|\mathbf h|^p,\\
 \norm{M_{W,p}^{\mathcal F}\mathbf h}_{L^p}^p
 &\le\sum_{Q\in\mathcal F}[W]_{A_p,Q}\int_Q|\mathbf h|^p\\
 &=\int_{\R^n}|\mathbf h(x)|^p\sum_{Q\in\mathcal F}[W]_{A_p,Q}\one_Q(x)\dd x\\
 &\le K[W]_{A_p}\norm{\mathbf h}_{L^p(\R^n,\C^d)}^p.
 \end{aligned}
\]
Combining this with \eqref{eq:knownstrong},
\[
 \norm{M_{W,p}^{\mathcal F}}_{L^p(\R^n,\C^d)\to L^p}
 \lesssim_{n,d,p}\min\left\{(K[W]_{A_p})^{1/p},[W]_{A_p}^{1/(p-1)}\right\}.
\]
By Chebyshev's inequality and the pointwise bound $M_{W,p}^{\mathcal F}\mathbf h\le M_{W,p}\mathbf h$,
\[
 \begin{aligned}
 \norm{M_{W,p}^{\mathcal F}}_{L^p(\R^n,\C^d)\to L^{p,\infty}}
 &\le\min\bigl\{\norm{M_{W,p}^{\mathcal F}}_{L^p(\R^n,\C^d)\to L^p},
 \norm{M_{W,p}}_{L^p(\R^n,\C^d)\to L^{p,\infty}}\bigr\}\\
 &\lesssim_{n,d,p}\min\{(K[W]_{A_p})^{1/p},
 [W]_{A_p}^{1/(p-1)},[W]_{A_p}^{2/p}\}.
 \end{aligned}
\]
Here the last term uses~\cite[Theorem~1.1]{LLOR} for $1<p<2$, and is redundant for $p\ge2$.
Taking the suprema over $[W]_{A_p}\le t$ and families satisfying \eqref{eq:depth} proves the upper estimates in \eqref{eq:depthstrong} and \eqref{eq:depthsharp}.
\end{proof}

\subsection{Sharp lower bounds when \texorpdfstring{$p\ge2$}{p >= 2}}\label{subsec:large}
\begin{proof}[Lower bound in \Thmref{thm:scalar}]
Let $t\ge K_{n,p}$ and choose $N=\lfloor(t/K_{n,p})^{1/(p-1)}\rfloor$ in \Lemref{lem:smooth}. Then $[w_*]_{A_p}\le t$ and
\[
 \norm{M_{w_*,p}^{\mathcal D}}_{L^p\to L^{p,\infty}}
 \ge \norm{M_{w_*,p}^{\mathcal D}f_*}_{L^{p,\infty}}\gtrsim_{n,p} N\eqsim_{n,p} t^{1/(p-1)}.
\]
For $1\le t<K_{n,p}$, take $w=1$ and a fixed nonnegative $L^p$ unit function $f_0\in C_c^\infty((0,1/2)^n)$. With $Q_0=[0,1/2)^n$,
\[
 \norm{M_{1,p}^{\mathcal D}f_0}_{L^{p,\infty}}
 \ge |Q_0|^{1/p}\langle f_0\rangle_{Q_0}\gtrsim_{n,p} t^{1/(p-1)}.
\]
Finally, \eqref{eq:conj} transfers the scalar lower bound to matrices, proving \eqref{eq:matrixmain}.
\end{proof}

\begin{proof}[Lower bounds in \Thmref{thm:depth} for $p\ge2$]
Fix $B=2^n$, $N\in\mathbb N$ and $v=w|_{\R_+}$ from \Propref{prop:Ap}, together with its unit function $h_C$. For each $K\ge1$ take $J=\min\{K,N\}$ and the chain
\[
 \mathcal I_K=\{I(Bs_j):N\le j<N+J\}.
\]
Since $B\ge2$, each $I(Bs_j)$ with $j\ge N$ contains $C$ and $E_j$. On the $J$ disjoint targets, \eqref{eq:sourceaverage} and \eqref{eq:testnorm} give
\[
 M_{v,p}^{\{I(Bs_j)\}}h_C(x)\ge B^{-1}S_C^{1/p'}\quad(x\in E_j),\qquad
 \left|\bigcup_{j=N}^{N+J-1}E_j\right|=JN.
\]
Consequently, by \eqref{eq:chainnorms},
\begin{equation}\label{eq:radlower}
 \norm{M_{v,p}^{\mathcal I_K}h_C}_{L^{p,\infty}}^p\ge B^{-p}S_C^{p-1}JN\gtrsim_{n,p}N^{p-1}\min\{K,N\}.
\end{equation}
By \eqref{eq:radlower} and \Lemref{lem:radial}, the lifted unit function $H(x)=h_C(q(x))\one_{\R_+^n}(x)$ and the weight $V(x)=v(q(x))$ satisfy
\[
 \begin{gathered}
 \mathcal F_K=\{Q(Bs_j):N\le j<N+J\},\qquad \#\mathcal F_K=J,\\
 [V]_{A_p}\le A_{n,p}N^{p-1},\qquad
 \norm{M_{V,p}^{\mathcal F_K}H}_{L^{p,\infty}}^p\gtrsim_{n,p}N^{p-1}J.
 \end{gathered}
\]
Here $A_{n,p}\ge1$ is independent of $N,K$. For $t\ge A_{n,p}$ take $N=\lfloor(t/A_{n,p})^{1/(p-1)}\rfloor$. Then $[V]_{A_p}\le t$ and
\[
 N^{p-1}\min\{K,N\}\eqsim_{n,p}t\min\{K,t^{1/(p-1)}\}.
\]
Consequently,
\[
 \begin{aligned}
 \norm{M_{V,p}^{\mathcal F_K}H}_{L^p}
 &\ge\norm{M_{V,p}^{\mathcal F_K}H}_{L^{p,\infty}}\\
 &\gtrsim_{n,p}\min\{(Kt)^{1/p},t^{1/(p-1)}\}.
 \end{aligned}
\]
This is the desired order for both norms when $p\ge2$. For $1\le t<A_{n,p}$, the constant weight $w=1$, unit function $\one_{[0,1)^n}$ and single cube $[0,1)^n$ give both norms equal to one. Embedding by $W=VI_d$ (or $I_d$ in the constant case) proves the matrix assertion.
\end{proof}

\subsection{Sharp lower bounds when \texorpdfstring{$1<p<2$}{1 < p < 2}}\label{sec:low}
\begin{proof}[Lower bounds in \Thmref{thm:depth} for $1<p<2$]
\emph{The weight and the common source.}
Fix $B=2^n$ and $m\in\mathbb N$, $m\ge m_0=\max\{100,8n\}$. Set
\[
 \varepsilon=m^{-1/(p-1)},\qquad \beta=(p-1)(1-\varepsilon),\qquad
 \tau=B^{-\lceil\log_B(1/\varepsilon)\rceil}.
\]
Then $\varepsilon/B<\tau\le\varepsilon\le1/m$. For $v$ from \Lemref{lem:patch}, its lift $V(x)=v(q(x))$ satisfies, by \Lemref{lem:radial},
\begin{equation}\label{eq:subcriticalbudget}
 [V]_{A_p(\R^n)}\le C_{n,p}m
\end{equation}
with $C_{n,p}\ge1$. Set $L_0=\lfloor\varepsilon^{-1}\rfloor$, $\delta=\tau B^{-L_0}$ and $C=(0,\delta)$. The source lies in all intervals used below. Its dual mass is
\[
 S=\int_C\sigma_v=\int_0^\delta(z/\tau)^{-1+\varepsilon}\dd z
 =\frac{\tau}{\varepsilon}B^{-L_0\varepsilon}.
\]
Since $L_0\varepsilon\le1$ and $B^{-1}<\tau/\varepsilon\le1$,
\[
 B^{-2}\le S\le1,\qquad \frac1{B\varepsilon}\le\frac S\tau\le\frac1\varepsilon,
 \qquad \left(\frac S\tau\right)^{p-1}\eqsim_{n,p}m.
\]
Use the normalized function $h_C$ of \Lemref{lem:common}.

\emph{The two parts of the chain.}
Choose $k_1<\cdots<k_{M_0}$ from
\[
 \{k\in\mathbb Z:\lceil m/2\rceil\le k\le m,\ n\mid(k+1)\},
 \qquad M_0\ge\left\lfloor\frac{m}{4n}\right\rfloor,
\]
where the count follows from the spacing $n$ and $m\ge8n$. By \eqref{eq:oldplateaux},
\[
 I_k^*\subseteq(0,2^{k+1}),\quad |I_k^*|=k\ge m/2,\quad v(z)^{1/p}=2^{k+1}\quad(z\in I_k^*).
\]
For the inner lengths and targets take
\[
 R_j=\delta B^j,\qquad A_j=(R_j/B,R_j)\quad(1\le j\le L_0).
\]
All averaging lengths lie in $B^{\mathbb Z}$ and exceed $\delta$. The intervals $A_j\subset(0,\tau)$ and $I_{k_i}^*\subset(1,\infty)$ are pairwise disjoint.

For each $K\ge1$ select at most $\lceil K/2\rceil$ outer and $\lfloor K/2\rfloor$ inner intervals. Set
\[
 J=\min\{\lceil K/2\rceil,\lfloor m/(4n)\rfloor\},\qquad
 \ell=\min\{\lfloor K/2\rfloor,L_0\}.
\]
The selected family is
\[
 \mathcal I_K=\{I(2^{k_i+1}):1\le i\le J\}\cup\{I(R_j):1\le j\le\ell\}.
\]
Thus
\[
 \mathcal I_K\subseteq\mathcal I_{K+1},\qquad \#\mathcal I_K=J+\ell\le K,
\]
with $v,h_C$ independent of $K$.

\emph{The outer intervals give the weak bound.}
By \eqref{eq:sourceaverage}, for $x\in I_{k_i}^*$,
\[
 M_{v,p}^{\mathcal I_K}h_C(x)\ge\frac{S^{1/p'}v(x)^{1/p}}{2^{k_i+1}}=S^{1/p'}.
\]
Hence, by \eqref{eq:chainnorms},
\[
 \norm{M_{v,p}^{\mathcal I_K}h_C}_{L^{p,\infty}}^p
 \ge S^{p-1}\sum_{i=1}^{J}|I_{k_i}^*|\gtrsim_{n,p}mJ .
\]

\emph{The inner intervals give the strong bound.}
For $1\le j\le\ell$ and $x\in A_j$, the smallest selected interval containing $x$ is $I(R_j)$. By \eqref{eq:chainidentity},
\[
 \begin{aligned}
 \int_{A_j}|M_{v,p}^{\mathcal I_K}h_C(x)|^p\dd x
 &= \frac{S^{p-1}}{R_j^p}\int_{R_j/B}^{R_j}(z/\tau)^\beta\dd z\\
 &=\frac{1-B^{-(1+\beta)}}{1+\beta}
   \left(\frac S\tau\right)^{p-1}
   \left(\frac{R_j}{\tau}\right)^{-(p-1)\varepsilon}
 \gtrsim_{n,p}\varepsilon^{-(p-1)}=m.
 \end{aligned}
\]
Since $R_j/\tau=B^{j-L_0}$,
\[
 1\le\left(\frac{R_j}{\tau}\right)^{-(p-1)\varepsilon}
 =B^{(L_0-j)(p-1)\varepsilon}\le B^{p-1}.
\]
Summing over $1\le j\le\ell$ gives
\[
 \norm{M_{v,p}^{\mathcal I_K}h_C}_{L^p}^p\gtrsim_{n,p}m\ell.
\]

\emph{Weak norm on the inner intervals.}
For $\mathcal P=\{I(R_j):1\le j\le L_0\}$, each selected interval containing $x$ has length at least $\max\{\delta,x\}$. By \eqref{eq:sourceaverage}, for $0<x<\tau$,
\[
 M_{v,p}^{\mathcal P}h_C(x)\le\frac{S^{1/p'}(x/\tau)^{\beta/p}}{\max\{\delta,x\}},
 \qquad M_{v,p}^{\mathcal P}h_C=0\quad\text{off }(0,\tau).
\]
Since $1+\beta-p=-(p-1)\varepsilon<0$,
\[
 \frac{x^{1+\beta}}{\max\{\delta,x\}^p}
 =\begin{cases}\delta^{-p}x^{1+\beta},&0<x\le\delta,\\x^{-(p-1)\varepsilon},&x>\delta\end{cases}
 \le\delta^{-(p-1)\varepsilon}.
\]
Consequently,
\[
 x|M_{v,p}^{\mathcal P}h_C(x)|^p\le (S/\tau)^{p-1}B^{L_0(p-1)\varepsilon}\lesssim_{n,p}m.
\]
For some $C_{n,p}'$ independent of $m$ and every $\lambda>0$,
\[
 \{x:M_{v,p}^{\mathcal P}h_C(x)>\lambda\}\subset(0,C_{n,p}'m/\lambda^p).
\]
Hence
\[
 \norm{M_{v,p}^{\mathcal P}h_C}_{L^{p,\infty}}^p\lesssim_{n,p}m,
 \qquad \norm{M_{v,p}^{\mathcal P}h_C}_{L^p}^p\gtrsim_{n,p}mL_0.
\]

\emph{Imposing the characteristic constraint.}
For $K\ge2$,
\[
 \begin{aligned}
 J&\ge\min\{K/2,m/(8n)\}\gtrsim_n\min\{K,m\},\\
 \ell&\ge\min\{K/3,\varepsilon^{-1}/2\}\gtrsim\min\{K,m^{1/(p-1)}\}.
 \end{aligned}
\]
For $K=1$, the outer interval gives both lower bounds by the domination of the weak norm by the strong norm. Thus, for all $K\ge1$,
\[
 \norm{M_{v,p}^{\mathcal I_K}h_C}_{L^{p,\infty}}^p\gtrsim_{n,p}m\min\{K,m\},\qquad
 \norm{M_{v,p}^{\mathcal I_K}h_C}_{L^p}^p\gtrsim_{n,p}m\min\{K,m^{1/(p-1)}\}.
\]
By \Lemref{lem:radial}, these inequalities transfer to dyadic cube chains of overlap $J+\ell\le K$, with the same lifted unit function and weight. For $t\ge2C_{n,p}m_0$, take $m=\lfloor t/C_{n,p}\rfloor\eqsim_{n,p}t$, so $[V]_{A_p}\le t$ by \eqref{eq:subcriticalbudget}. For smaller $t\ge1$, use $w=1$, $h=\one_{[0,1)^n}$ and $\mathcal F_K=\{[0,1)^n\}$. Both norms equal one. The embedding $W=VI_d$ and $h\mathbf e$ proves \Thmref{thm:depth}.
\end{proof}

\begin{proof}[Proof of \Corref{cor:ratio}]
Take $t=[W]_{A_p}\ge1$ and use \eqref{eq:depth} in \Thmref{thm:depth}. For the common examples, $K=1$ gives
\[
 t\lesssim_{n,p}\norm{M_{w_t,p}^{\mathcal F_1}h_t}_{L^p}^p\le[w_t]_{A_p}\le t.
\]
Thus $[w_t]_{A_p}\eqsim_{n,p}t$. The lower bounds on $h_t$ and the upper bounds in \Thmref{thm:depth} give
\[
 \begin{aligned}
 \norm{M_{w_t,p}^{\mathcal F_K}h_t}_{L^p}^p
 &\eqsim_{n,p}t\min\{K,t^{1/(p-1)}\},\\
 \norm{M_{w_t,p}^{\mathcal F_K}h_t}_{L^{p,\infty}}^p
 &\eqsim_{n,p}t\min\{K,t^{1/(p-1)},t\}.
 \end{aligned}
\]
The corresponding operator norms have the same orders. Substituting $[w_t]_{A_p}\eqsim_{n,p}t$ proves the sharpness of \eqref{eq:actualstrong} and \eqref{eq:actualweak}. Dividing the two expressions gives \eqref{eq:ratio}.
\end{proof}

\section{Proofs of the auxiliary lemmas}\label{sec:aux}
\subsection{The common test and single averages}
\begin{proof}[Proof of \Lemref{lem:common}]
By the definition of $h_C$ and \eqref{eq:scalarop},
\[
 \int_C|h_C|^p=S^{-1}\int_C\sigma_v=1,\qquad \int_C v^{-1/p}h_C=S^{1/p'}.
\]
H\"older's inequality gives
\[
 \int_Cv^{-1/p}|f|\le S^{1/p'}\norm{f}_{L^p},
\]
with equality for $h_C$. For $I\supset C$ up to a null set, this proves \eqref{eq:sourceaverage}. On $r_{j-1}<x<r_j$ the smallest eligible interval is $I(r_j)$. Every eligible interval contains $C$, so
\[
 M_{v,p}^{\mathcal I}h_C(x)=\max_{j\le k\le J}\frac{S^{1/p'}v(x)^{1/p}}{r_k}
 =\frac{S^{1/p'}v(x)^{1/p}}{r_j}.
\]
Integrating the $p$th power proves \eqref{eq:chainidentity}. For $a>0$,
\[
 \bigl|\{M_{v,p}^{\mathcal I}h_C>aS^{1/p'}\}\bigr|
 =\sum_{j=1}^J\bigl|\{x\in(r_{j-1},r_j):v(x)>(ar_j)^p\}\bigr|.
\]
Taking the supremum gives \eqref{eq:weakchain}. On each $E_j$, $M_{v,p}^{\mathcal I}h_C\ge S^{1/p'}v^{1/p}/r_j$. Therefore disjointness gives \eqref{eq:strongsource}. Under the additional hypothesis in \eqref{eq:chainnorms}, for every $0<\lambda<aS^{1/p'}$,
\[
 \lambda^p\bigl|\{M_{v,p}^{\mathcal I}h_C>\lambda\}\bigr|
 \ge\lambda^p\sum_{j=1}^J|E_j|.
\]
Letting $\lambda\uparrow aS^{1/p'}$ proves \eqref{eq:chainnorms}.
\end{proof}
For a scalar $A_p$ weight $v$, the H\"older equality on a cube $Q$ uses $\sigma_v^{1/p}\one_Q/(\int_Q\sigma_v)^{1/p}$. Then
\begin{align*}
 \norm{M_{v,p}^{\{Q\}}}_{L^p\to L^p}^{p}&=\langle v\rangle_Q\langle\sigma_v\rangle_Q^{p-1},\\
 \norm{M_{v,p}^{\{Q\}}}_{L^p\to L^{p,\infty}}^{p}&=\frac{\norm{v\one_Q}_{L^{1,\infty}}}{|Q|}\langle\sigma_v\rangle_Q^{p-1}.
\end{align*}
Both products are at most $[v]_{A_p}$. The supremum over $Q$ of the second product is $[v]_{A_p^*}$ in Sweeting~\cite[Definition~2.1]{Sweeting}. Its finiteness characterizes scalar multiplier weak boundedness by~\cite[Theorem~1.2]{Sweeting}.

\subsection{Reduction of the strong matrix bound to dimension one}
\begin{proof}[Proof of \Lemref{lem:dim}]
Fix $\mathbf f\in L^p(\R,\C^d)$ and $R>0$. On $\R^2$ set $\widehat W(x,y)=W(x)$ and $F_R(x,y)=\mathbf f(x)\one_{(0,R)}(y)$. By Fubini,
\[
 [\widehat W]_{A_p(\R^2)}=[W]_{A_p(\R)},\qquad
 \norm{F_R}_{L^p(\R^2,\C^d)}^p=R\norm{\mathbf f}_{L^p(\R,\C^d)}^p.
\]
Let $M_{W,p}^{\le R}\mathbf f(x)=\sup_{I\ni x,\,|I|\le R}M_{W,p}^{\{I\}}\mathbf f(x)$. If $x\in I$, $|I|\le R$ and $0<y<R$, choose an interval $J\subseteq[0,R]$ containing $y$ with $|J|=|I|$. The cube $I\times J$ gives
\[
 M_{\widehat W,p}^{\{I\times J\}}F_R(x,y)=M_{W,p}^{\{I\}}\mathbf f(x),\qquad
 M_{\widehat W,p}F_R(x,y)\ge M_{W,p}^{\le R}\mathbf f(x).
\]
Integrating over the strip $\R\times(0,R)$ and using the planar bound, with constant $C_{2,d,p}$, we obtain
\[
 R\norm{M_{W,p}^{\le R}\mathbf f}_{L^p(\R)}^p
 \le\norm{M_{\widehat W,p}F_R}_{L^p(\R^2)}^p
 \le C_{2,d,p}^p[W]_{A_p(\R)}^{p/(p-1)}R\norm{\mathbf f}_{L^p(\R,\C^d)}^p.
\]
Divide by $R$. Since $M_{W,p}^{\le R}\mathbf f\uparrow M_{W,p}\mathbf f$ almost everywhere as $R\to\infty$, monotone convergence gives the claimed estimate.
\end{proof}

\subsection{Lifting to cubes}\label{subsec:liftproof}
\begin{proof}[Proof of \Lemref{lem:radial}]
Fix an axis-parallel cube $Q$ of side length $\ell$. Write $r(x)=|x|_\infty$,
$r_0=\inf_Qr$, and $r_1=\sup_Qr$. As $r$ is $1$-Lipschitz in the supremum norm, $r_1-r_0\le\ell$. We find an interval $J$ such that
\begin{equation}\label{eq:radavg}
 \frac{1}{|Q|}\int_Q h(r(x)^n)\dd x\le c_n\langle h\rangle_J
 \qquad (h\ge0\text{ measurable}).
\end{equation}
When $r_0\le\ell$, take $J=(0,r_1^n)$. The cube $Q$ lies in
$[-r_1,r_1]^n$, and $|\{x:r(x)^n\le z\}|=2^nz$. Consequently
\[
 \frac{1}{|Q|}\int_Q h(r(x)^n)\dd x\le
 \frac{2^n}{\ell^n}\int_0^{r_1^n}h(z)\dd z
 \le 4^n\langle h\rangle_J ,
\]
where the last inequality uses $r_1\le2\ell$.

When $r_0>\ell$, take $J=(r_0^n,r_1^n)$, which has positive length. Partition $Q$ up to null sets by a maximizing coordinate and sign. Integrate first in $x_i$, substituting $z=|x_i|^n$ with Jacobian at most $(nr_0^{n-1})^{-1}$. The transverse box has volume $\ell^{n-1}$. Summing over $2n$ pieces gives
\[
 \int_Q h(r(x)^n)\dd x
 \le\frac{2\ell^{n-1}}{r_0^{n-1}}\int_J h(z)\dd z .
\]
Since $r_1\le r_0+\ell<2r_0$,
\[
 |J|=r_1^n-r_0^n\le n(2r_0)^{n-1}\ell,
 \qquad \frac{1}{\ell^n}\int_Qh(r(x)^n)\dd x\le n2^n\langle h\rangle_J.
\]
Thus \eqref{eq:radavg} holds with $c_n=\max\{4^n,n2^n\}$. Applying it to $v$ and $v^{-1/(p-1)}$ on the same $J$ gives
\[
 (\langle V\rangle_Q)(\langle V^{-1/(p-1)}\rangle_Q)^{p-1}
 \le c_n^p(\langle v\rangle_J)(\langle v^{-1/(p-1)}\rangle_J)^{p-1}\le c_n^p[v]_{A_p(\R_+)}.
\]
It remains to prove the norm identities. Since $|\{x\in\R_+^n:q(x)<z\}|=z$, for nonnegative measurable $G$ we have
\begin{equation}\label{eq:pushforward}
 \int_{\R_+^n}G(q(x))\dd x=\int_0^\infty G(z)\dd z.
\end{equation}
This gives the identity for the unweighted norm and, for each $r_j$,
\[
 \frac{1}{r_j}\int_{Q(r_j)}V(y)^{-1/p}|F(y)|\dd y=\frac{1}{r_j}\int_0^{r_j}v(z)^{-1/p}|f(z)|\dd z.
\]
On $\R_+^n$, $x\in Q(r_j)$ exactly when $q(x)<r_j$. Hence
\[
 M_{V,p}^{\mathcal F}F(x)=M_{v,p}^{\mathcal I}f(q(x))\quad\text{a.e. on }\R_+^n.
\]
The function on the left vanishes off the positive orthant, up to a null set. Hence \eqref{eq:pushforward} gives
\[
 \begin{aligned}
 \int_{\R^n}|M_{V,p}^{\mathcal F}F|^p&=\int_0^\infty|M_{v,p}^{\mathcal I}f|^p,\\
 |\{M_{V,p}^{\mathcal F}F>\lambda\}|&=|\{M_{v,p}^{\mathcal I}f>\lambda\}|\quad(\lambda>0).
 \end{aligned}
\]
This proves \eqref{eq:transfer}. Since $r_j\in B^{\mathbb Z}$, $Q(r_j)\in\mathcal D$. For a finite chain of $J$ cubes, the smallest has positive volume and the overlap is exactly $J$.
\end{proof}

\subsection{The characteristic of the periodic weight}\label{sec:periodic-characteristic}
\begin{proof}[Proof of \Propref{prop:Ap}]
With the intervals in \eqref{eq:CE}, write $\sigma=w^{-1/(p-1)}$ and $\Phi(z)=(1+\operatorname{dist}(z,\Gamma))^{-1}$ on $[0,L]$. We have $256BN\le b<256B^2N$ and $\log(1+L)\lesssim_B N$. Continuity, positivity and periodicity follow from \eqref{eq:weight}.

\emph{The source and target estimates.} If $z\in C_k$, then $5s_k/4\le z\le3s_k/2$. The gap in $\Gamma$ containing $C_k$ has endpoints $s_k+N$ and $Bs_k$. Thus $z-s_k-N\ge s_k/4-N\ge s_k/5$ and $Bs_k-z\ge s_k/2$. Since $r(z)=z$,
\[
 Ns_k^{p-1}\le w(z)\le B^p(N+5)s_k^{p-1},\qquad
 \int_{C_k}\sigma\eqsim_{B,p}|C_k|\,(Ns_k^{p-1})^{-1/(p-1)}\eqsim_{B,p}N^{-1/(p-1)}.
\]
Consequently,
\[
 S_C=\sum_{k=0}^{N-1}\int_{C_k}\sigma\eqsim_{B,p}N^{1-1/(p-1)},
 \qquad w|_{E_j}\ge s_j^p,\qquad |E|=\sum_{j=N}^{2N-1}|E_j|=N^2.
\]
This proves \eqref{eq:testnorm}. Since $C,E_{2N-1}\subset[0,L]$ and $L=Bs_{2N-1}$,
\[
 \langle w\rangle_{[0,L]}\langle\sigma\rangle_{[0,L]}^{p-1}\ge L^{-p}Ns_{2N-1}^p S_C^{p-1}\gtrsim_{B,p}N^{p-1}.
\]

\emph{Short intervals away from zero.} Let $H=[u,u+\ell]\subseteq[0,L]$ with $u\ge8\ell>0$. Since $(u+\ell)/u\le9/8<B$, $H$ meets at most two consecutive pieces of the partition $[0,b),J_k=[s_k,Bs_k)$ ($0\le k<2N$), and $\sup_H r\le(9/8)\inf_H r$. Let $d_H=\inf_H\operatorname{dist}(\cdot,\Gamma)$. The distance is $1$-Lipschitz. If $\ell\le1$ or $d_H\ge\ell$, then
\[
 \sup_H\Phi/\inf_H\Phi\le(1+d_H+\ell)/(1+d_H)\le2.
\]
Otherwise, $\ell>1$ and $d_H<\ell$. On each gap $(s_k+N,Bs_k)$ the distance equals the minimum of the distances to its two endpoints. On each half-gap a subinterval of length at most $\ell$ has $\Phi$-integral at most $\log(1+\ell)$. Since $H$ meets at most two pieces, the initial interval contributes at most $b$ and the short targets at most $2N$. Hence
\[
 \int_H\Phi\le b+2N+4\log(1+\ell)\lesssim_B N,\qquad
 \frac{\langle\Phi\rangle_H}{\inf_H\Phi}\le \frac{1+2\ell}{\ell}\int_H\Phi\lesssim_B N.
\]
The latter uses $\sup_H\operatorname{dist}(\cdot,\Gamma)\le d_H+\ell<2\ell$. In both cases, writing $r_H=\inf_Hr$ and $\phi_H=\inf_H\Phi$, we get
\[
 \begin{aligned}
 \langle w\rangle_H&\lesssim_{B,p}Nr_H^{p-1}+Nr_H^p\phi_H\lesssim_{B,p}N\inf_H w,\\
 \langle w\rangle_H\langle\sigma\rangle_H^{p-1}&\le\langle w\rangle_H/\inf_H w\lesssim_{B,p}N\le N^{p-1}.
 \end{aligned}
\]

\emph{Initial intervals.} On $J_k$, with $s=s_k$, the gap calculation gives
\[
 \int_{J_k}\Phi=N+2\log(1+((B-1)s-N)/2)\lesssim_B N.
\]
Since $s\le r\le Bs$, \eqref{eq:weight} gives
\[
 \begin{aligned}
 \int_{J_k}w
 &\le N(Bs)^{p-1}(B-1)s+(Bs)^p\int_{J_k}\Phi
 \lesssim_{B,p}Ns^p,\\
 \int_{J_k}\sigma
 &\le (Ns^{p-1})^{-1/(p-1)}(B-1)s
 =(B-1)N^{-1/(p-1)}.
 \end{aligned}
\]
On $[0,b]$, $w=b^{p-1}(N+b)$, so
\[
 \int_0^b w=b^p(N+b)\eqsim_{B,p}Nb^p,
 \qquad \int_0^b\sigma=(N+b)^{-1/(p-1)}.
\]
For $b\le R\le L$, summing the geometric series and at most $2N$ dual masses gives
\begin{equation}\label{eq:mass}
 \begin{aligned}
 \int_0^Rw&\lesssim_{B,p}Nb^p+N\!\sum_{s_k\le R}s_k^p\lesssim_{B,p}NR^p,\\
 \int_0^R\sigma&\le(N+b)^{-1/(p-1)}+2N(B-1)N^{-1/(p-1)}\lesssim_B N^{(p-2)/(p-1)}.
 \end{aligned}
\end{equation}
For $H=[u,u+\ell]\subseteq[0,L]$ with $u<8\ell$, take $R=\min\{9\ell,L\}$. Then $H\subseteq[0,R]$. The weight is constant if $R<b$. Otherwise, by \eqref{eq:mass},
\[
 \langle w\rangle_H\langle\sigma\rangle_H^{p-1}\le\ell^{-p}\int_0^Rw\left(\int_0^R\sigma\right)^{p-1}\lesssim_{B,p}N^{p-1}(R/\ell)^p\lesssim_{B,p}N^{p-1}.
\]
Thus all subintervals of $[0,L]$ have the required bound.

\emph{Reflection and periodicity.} Set $\rho(x)=\operatorname{dist}(x,2L\mathbb Z)$. On consecutive intervals between points of $L\mathbb Z$, $\rho$ is an isometry onto $[0,L]$ up to endpoints. For an interval $H$ with $|H|<L$, the image $J=\rho(H)$ is an interval of length at most $|H|$ with multiplicity at most two. For $|H|\ge L$, enlarge $H$ to endpoints in $L\mathbb Z$ at cost at most $3|H|$ and take $J=[0,L]$. In either case,
\[
 \langle h\circ\rho\rangle_H\le3\langle h\rangle_J\qquad(h\ge0).
\]
Applying the preceding estimate to $w_0$ and $w_0^{-1/(p-1)}$ with the same $J$ yields
\[
 \langle w\rangle_H\langle\sigma\rangle_H^{p-1}
 \le3^p\langle w_0\rangle_J\langle w_0^{-1/(p-1)}\rangle_J^{p-1}
 \lesssim_{B,p}N^{p-1}.
\]
The lower bound above completes \eqref{eq:Apglobal}.
\end{proof}

\subsection{Smooth periodic weights in arbitrary dimensions}
\begin{proof}[Proof of \Lemref{lem:smooth}]
Use \Propref{prop:Ap} with $B=2$ and $h_C=S_C^{-1/p}\sigma_w^{1/p}\one_C$. For $N\le j<2N$, the dyadic interval $[0,2s_j)$ contains $C$ and $E_j$. By \eqref{eq:sourceaverage} and \eqref{eq:testnorm}, with $\lambda=\tfrac12S_C^{1/p'}$,
\[
 \begin{gathered}
 M_{w,p}^{\mathcal D}h_C(x)\ge\frac{S_C^{1/p'}w(x)^{1/p}}{2s_j}\ge\lambda\quad(x\in E_j),\\
 \lambda|E|^{1/p}\eqsim_p N^{(p-2)/p}N^{2/p}=N.
 \end{gathered}
\]
Set $w_n(x)=w(x_1)$, so $[w_n]_{A_p}=[w]_{A_p(\R)}$ by Fubini. For $n\ge2$ take
\[
 H(x)=L^{-(n-1)/p}h_C(x_1)\one_{[0,L)^{n-1}}(x'),\qquad x=(x_1,x').
\]
For $n=1$ use $H=h_C$, $E_1=E$ and $\lambda_1=\lambda$. Fubini gives $\norm H_{L^p}=1$. For $n\ge2$, partition $[0,L)^{n-1}$ into dyadic cubes of side $2s_j$. On each product with $[0,2s_j)$ the average of $w_n^{-1/p}H$ is $L^{-(n-1)/p}S_C^{1/p'}/(2s_j)$. Set $E_n=E\times[0,L)^{n-1}$ for $n\ge2$. Then
\[
  M_{w_n,p}^{\mathcal D}H\ge\lambda_n=L^{-(n-1)/p}\lambda\quad\text{on }E_n,
 \qquad \lambda_n |E_n|^{1/p}=\lambda |E|^{1/p}.
\]

Choose a nonnegative $\eta\in C_c^\infty((-1,1))$ with $\int\eta=1$, $\eta_\epsilon(x)=\epsilon^{-1}\eta(x/\epsilon)$, and $u(x)=\int\eta_\epsilon(y)w(x-y)\dd y$. By continuity and the positive minimum of periodic $w$, $w/2\le u\le2w$ for sufficiently small $\epsilon>0$. Hence $U(x)=u(x_1)$ satisfies $[U]_{A_p}\le4[w]_{A_p}\le K_{n,p} N^{p-1}$ for some $K_{n,p}\ge1$. Since $w_n/2\le U\le2w_n$, the multiplier ratios satisfy
\[
 U(x)^{1/p}U(y)^{-1/p}\ge4^{-1/p}w_n(x)^{1/p}w_n(y)^{-1/p}.
\]
Thus $M_{U,p}^{\mathcal D}H\ge4^{-1/p}\lambda_n$ on $E_n$, with the same $L^p$ unit function $H$. Choose nonnegative $L^p$ unit functions $G_j\in C_c^\infty((0,L)^n)$ with $G_j\to H$ in $L^p$. For every cube $Q$ used to obtain the lower bound on $E_n$,
\[
 \left|\int_Q U^{-1/p}(G_j-H)\right|
 \le\norm{U^{-1/p}\one_Q}_{L^{p'}}\norm{G_j-H}_{L^p}\longrightarrow0.
\]
There are only finitely many such cubes, and $\int_Q U^{-1/p}H>0$ for each. Thus one approximant $G$ satisfies
\[
 \int_Q U^{-1/p}G\ge\frac12\int_Q U^{-1/p}H\quad\text{for every such }Q,
 \qquad M_{U,p}^{\mathcal D}G\ge\frac{\lambda_n}{2\,4^{1/p}}\quad\text{on }E_n.
\]

Finally, $P=2L$ is a power of two. Set $\widetilde w(x)=u(Px_1)$ and $G_P(x)=P^{n/p}G(Px)$. Since $P$ is a dyadic scale,
\[
 [\widetilde w]_{A_p}=[U]_{A_p},\quad \norm{G_P}_{L^p}=1,\quad
 M_{\widetilde w,p}^{\mathcal D}G_P(x)=P^{n/p}M_{U,p}^{\mathcal D}G(Px).
\]
Hence $\norm{M_{\widetilde w,p}^{\mathcal D}G_P}_{L^{p,\infty}}=\norm{M_{U,p}^{\mathcal D}G}_{L^{p,\infty}}\ge 2^{-1}4^{-1/p}\lambda_n|E_n|^{1/p}\gtrsim_{n,p} N$. Take $u_*(z)=u(Pz)$, $w_*=\widetilde w$ and $f_*=G_P$.
\end{proof}

\subsection{The characteristic after the inner modification}
\begin{proof}[Proof of \Lemref{lem:patch}]
Write $\sigma_v=v^{-1/(p-1)}$. Since $u_m=1$ on $[0,1]$ and $\tau\le\varepsilon\le1/m$, we have $\int_0^\tau\sigma_v=\tau/\varepsilon\le1$ and $\tau<1/4$.
For $0<b\le1$, the $A_p$ product on $(0,b)$ satisfies
\[
 (\langle v\rangle_{(0,b)})(\langle\sigma_v\rangle_{(0,b)})^{p-1}
 \le\begin{cases}
 (1+\beta)^{-1}\varepsilon^{-(p-1)},&b\le\tau,\\
 \bigl((\tau/\varepsilon+b-\tau)/b\bigr)^{p-1},&\tau<b\le1.
 \end{cases}
\]
The first line is exact, and the second is bounded by $\varepsilon^{-(p-1)}$ because
\[
 \frac{\tau/\varepsilon+b-\tau}{b}
 =1+\frac{\tau}{b}(\varepsilon^{-1}-1)\le\varepsilon^{-1}.
\]
Hence both bounds are at most $m$. If $0\le a<b\le1$ and $a\le b/2$, comparison with $(0,b)$ gives
\[
 \langle v\rangle_{(a,b)}\langle\sigma_v\rangle_{(a,b)}^{p-1}
 \le2^p\langle v\rangle_{(0,b)}\langle\sigma_v\rangle_{(0,b)}^{p-1}\le2^pm.
\]
If $b/2<a<b\le1$, then $v(z)=\min\{z/\tau,1\}^{\beta}$ on $(a,b)$, and
\[
 \langle v\rangle_{(a,b)}\langle\sigma_v\rangle_{(a,b)}^{p-1}
 \le\frac{\sup_{(a,b)}v}{\inf_{(a,b)}v}\le(b/a)^\beta\le2^\beta.
\]
If $a\ge\tau$, the weight agrees with $u_m$ on $(a,b)$. It remains to consider $a<\tau$ and $b>1$. Since $v\le u_m$ and $\tau<1/4$,
\[
 \int_a^b u_m^{-1/(p-1)}\ge1-\tau\ge3/4,
 \qquad \int_a^b\sigma_v\le\int_a^b u_m^{-1/(p-1)}+\tau/\varepsilon
 \le\frac73\int_a^b u_m^{-1/(p-1)}.
\]
Consequently,
\[
 \langle v\rangle_{(a,b)}\langle\sigma_v\rangle_{(a,b)}^{p-1}
 \le(7/3)^{p-1}[u_m]_{A_p}\lesssim_p m.
\]
Taking the supremum over all intervals in $\R_+$ completes the proof.
\end{proof}

\medskip
\noindent{\bf Funding }
Not applicable.
\medskip

\noindent{\bf Data Availability} This manuscript has no associated data.

\medskip
\noindent{\bf\Large Declarations}
\medskip

\noindent\textbf{Conflicts of Interest}\space The authors declare no conflicts of interest.

\medskip
\noindent{\bf AI disclosure}
During the preparation of this manuscript, the authors used ChatGPT and Codex to discuss possible proof strategies, assist in drafting and checking parts of the mathematical arguments, organize references, and improve the clarity and readability of the exposition. The authors take full responsibility for the mathematical statements, proofs, and exposition in this manuscript.

\end{document}